\documentclass[journal, two-sided]{IEEEtran}
\usepackage{algorithmic}

\usepackage{algorithm}
\usepackage{mathrsfs,amsthm}
\IEEEoverridecommandlockouts
\makeatletter
\def\ps@headings{%
\def\@oddhead{\mbox{}\scriptsize\rightmark \hfil \thepage}%
\def\@evenhead{\scriptsize\thepage \hfil \leftmark\mbox{}}%
\def\@oddfoot{}%
\def\@evenfoot{}}
\makeatother 
\usepackage{graphicx,amssymb,stfloats,subfigure,multicol}
\usepackage{epsfig,epsf,psfrag,amssymb,amsfonts,latexsym,graphicx,mathrsfs,subfigure}
\usepackage{bbm}
\usepackage{chngpage}
\usepackage[dvips]{color}
\usepackage{textcomp}
\usepackage{hyperref}
\usepackage[hyphenbreaks]{breakurl}
\usepackage[normalem]{ulem}
\newtheorem{Def}{Definition}

\newtheorem{thm}{Theorem}
\newtheorem{lem}{Lemma}
\newtheorem{prop}{Proposition}

 \def\old#1{}    

\def\beq{\begin{equation}}
\def\eeq{\end{equation}}
\def\bea{\begin{eqnarray}}
\def\eea{\end{eqnarray}}
\def\ba{\begin{array}}
\def\ea{\end{array}}
\def\bitem{\begin{itemize}}
\def\eitem{\end{itemize}}
\def\benum{\begin{enumerate}}
\def\eenum{\end{enumerate}}
\def\defeq{{\stackrel{\Delta}{=}}}
\def\edoc{\end{document}}

\def\ie{{\it i.e.,\ \/}}

\definecolor{bgrd}{rgb}{1,1,1}
\definecolor{grey}{rgb}{0.9,0.9,0.6}
\definecolor{gray}{rgb}{0.5,0.5,0.5}
\definecolor{dkr}{rgb}{0.6,0.2,0.2}
\definecolor{dkg}{rgb}{0,0.5,0}
\definecolor{dkb}{rgb}{0.0,0.1,0.7}
\definecolor{light-gray}{gray}{0.85}

\begin{document}

\title{Large-scale Charging of Electric Vehicles:\\
A Multi-Armed Bandit Approach}
\author{\Large Zhe Yu$^\dagger$, Yunjian Xu$^\ddagger$, and Lang Tong$^\dagger$
\thanks{\scriptsize
Z. Yu$^\dagger$ and L. Tong$^\dagger$ are with the School of Electrical and Computer
Engineering, Cornell University, Ithaca, NY 14853, USA. Y. Xu$^\ddagger$ is with the School of Engineering Systems and Design Pillar, Singapore University of Technology and Design, Singapore, 487372. Email:
{\tt\{zy73,lt35\}@cornell.edu, yunjian\_xu@sutd.edu.sg}. This work is supported in part by the National Science Foundation under
Grant CNS-1248079. Part of this work is presented at the 2015 Allerton conference in October 2015.}}

\maketitle

\begin{abstract}
The successful launch of electric vehicles (EVs) depends critically on the availability of convenient and economic charging facilities. The problem of scheduling of large-scale charging of EVs by a service provider is considered.   A Markov decision process model is introduced in which EVs arrive randomly at a charging facility with random demand and completion deadlines.  The service provider faces random charging costs, convex non-completion penalties, and a  peak power constraint that limits the maximum  number of simultaneous activation of EV chargers.

 Formulated as a restless multi-armed bandit problem, the EV charging problem is shown to be indexable.  A closed-form expression of the Whittle's index is obtained for the case when the charging costs are constant.  The Whittle's index policy, however, is not optimal in general.  An enhancement of the Whittle's index policy based on spatial interchange according to  the less laxity and longer processing time principle is presented.  The proposed policy outperforms existing charging algorithms, especially when the charging costs are time varying.
\end{abstract}

\begin{IEEEkeywords}
Charging of electric vehicles; deadline scheduling; Markov decision processes; multi-armed bandit problem;
 Whittle's index.
\end{IEEEkeywords}

\section{Introduction}


\IEEEPARstart{A}{ccording} to a recent study \cite{Creutzig&Etal:S2015}, a transition from gasoline based transportation to electric vehicles (EV), coupled with integrating renewable resources for EV charging, will play a critical role in achieving the goal of halving the $\mbox{CO}_2$ emissions by 2050. In 2015, the global electric drive vehicle sales (including hybrid, plug-in vehicle and battery electric vehicle) reaches over 498,000, which is $2.87\%$ of the annually all vehicle sales \cite{EDTA:15}. In the US, the annually EV sales has grown 20 times since 2011. Similar trend exists in the EV charging station market. Through the end of 2014, there are more than 15,000 fast charging points and 94,000 slow charging points over the world. The EV charging station stock more than doubled for slow charging points between the end of 2012 and 2014, and increased eightfold for fast charging points \cite{IEA:15}. As of March 2016, there are more than 12,700 electric stations and 31,800 charging outlets deployed in the United States \cite{DoE:15}.

EV charging services play an essential role in the successful launch of EVs. A sufficient amount of charging services attracts more consumers to purchase EVs and high EV market share brings more investment in the charging services \cite{Yu&Li&Tong:14Allerton,Winebrake&Farrell:TRD97}.
Large charging facilities with fast charging capabilities in public spaces such as parking garages, parking lots at commercial locations, and highway rest stops serve to alleviate range anxiety of EV consumers and stimulate the market share of EVs. These facilities that serve a large number of EVs at any given time bring the additional benefit of providing ancillary services and maintaining operation stability of the power grid \cite{Kempton&Tomic:05JPS, Brooks:book}.

Large scale EV charging at the capacity of hundreds of vehicles faces a different set of technical  challenges from those associated with individual home charging.  First,  consumers expect charging to be completed within a relatively short period of time.  Thus, fast charging devices operated at high peak power becomes essential. Currently, level 2 and DC fast charging are most widely used in public charging stations. Level 2 charging supplies up to 30 miles of travel for one hour of charging with a 6.6kWh on-board charger. DC fast charging supplies up to 40 miles of range of driving for every 10 minutes of charging, which equals approximately 15 average size residential central air conditioning units.
 These types of charging, if un-managed, may have detrimental effects on power system reliability \cite{Clement&Etal:TPS2010,Sortomme&Etal:TSG2011}.  It is thus necessary to limit the number of simultaneously activated chargers.

Second, there is a high level of uncertainty in charging demand at public facilities. EVs arrive at a charging facility randomly, each with stochastic demand and random deadlines, which makes it difficult for the scheduler to meet  consumer demands.

Third, the cost (or the profit) of the service provider may be stochastic.  For instance, the service provider may participate in the wholesale electricity market and is subject to real-time  price fluctuations.  In addition, the service provider may integrate local renewable energy such as solar with intermittent generatioin.

Finally, the energy management system that schedules EV charging operates in real time. Therefore, the scheduling algorithm must be scalable with respect to the size of the charging facility, which rules out the use of brute-force optimization techniques.

\subsection{Summary of Results}
This work extends the results from \cite{Yu&Xu&Tong:2015Allerton}, which is the first article to apply the restless multi-armed bandit problems to the EV charging.
We propose an online  scheduling algorithm that is computationally scalable and capable of dealing with demand and cost uncertainty.   We introduce a constrained Markov decision process (MDP) model with the objective of  maximizing expected (discounted) profit subject to a constraint on the maximum number of simultaneously activated  chargers.
The constructed MDP model captures the randomness in EV arrivals, EVs' charging requests and deadlines, as well as the charging costs.
The evolution of charging cost is random, and is assumed to be independent of the actions taken by
the operator.

 We note that computing exact optimal scheduling policies by brute-force dynamic programming is
intractable, because the number of system states grows exponentially with the number of chargers.
In order to derive effective online  scheduling algorithms,
we reformulate the MDP as a restless multi-armed bandit (MAB) problem with simultaneous  plays \cite{Whittle:1988JAP}.
We first establish the indexability of the formulated restless MAB problem, which enables us to apply
the {\it Whittle's index policy} to the EV charging problem.
The special structure of the EV charging problem, in particular the pre-determined charging deadlines, simplifies the computation of the Whittle's index.  For the case with constant charging cost we derive the Whittle's indexes in closed form.

We establish the optimality of the Whittle's index policy for random charging cost when the constraint on the number of simultaneously activated EV charging is loose. When the constraint of simultaneous activation is  strict, Whittle's index policy is not optimal in general. In this context, we provide a procedure based on the LLLP (Less Laxity and Longer Processing) principle [11] as an improvement of the Whittle's index policy.
  Numerical results demonstrate that
the LLLP principle could significantly improve the performance of the  Whittle's index policy,
   especially when the charging cost is stochastic and the EV arrival traffic is relatively heavy.

\subsection{Related Work}

The centralized EV charging problem considered in this paper falls in the category of multi-processor deadline scheduling problems.  In this context, EVs are jobs and chargers are processors.

Earlier work on deadline scheduling are based on the deterministic worst case  objectives.  The problem of deadline scheduling with one processor is well understood.  In this case, simple online algorithms such as the earliest deadline first (EDF) policy \cite{Liu&Layland:73ACM} and the least-laxity first
(LLF) policy \cite{Dertouzos:74IFIPC}, are optimal, when the completion of all tasks  before deadlines
is feasible. Under certain conditions, it is shown that the EDF scheduling minimizes the amount of unfinished work in single-processor deadline scheduling \cite{Panwar&Towsley&Wolf:88JACM,Towsley&Panwar:90RT}.
 There is also a substantial literature on
deadline scheduling of multiple processors
(for a survey, see \cite{DB11}). It is shown in \cite{DerTouzos&Mok:89TSE}
that an optimal online scheduling policy does not exist in general for the worst case performance measure.

The problem of stochastic multi-processor deadline scheduling, of which the EV charging problem is a special case, is less understood, primarily because the stochastic dynamic programming for such problems are not tractable in practice.  The work similar to ours are \cite{BE89,BT97} where the authors studied the deadline scheduling problem in wireless communications. The authors of \cite{Raghunathan&Etal:2008INFOCOM} formulate the stochastic deadline scheduling problem (in wireless communications) as a restless MAB problem,  and establish indexability for the formulated MAB problem.
      Related problems of scheduling packets with deadlines in ad hoc networks are studied in
       \cite{Jaramillo&Srikant&Ying:2011SAC}.
       We note that there are fundamental differences between the job arrival and processing cost models
       adopted in this paper and in the aforementioned literature, and that the results derived in these existing works
       do not apply to our model.

       \old{
  There are several nontrivial differences between the models in \cite{Raghunathan&Etal:2008INFOCOM,Jaramillo&Srikant&Ying:2011SAC} and that in the current paper.  For instance, the arrival models used in  \cite{Raghunathan&Etal:2008INFOCOM}  are either simultaneous or periodic.  The cost models in \cite{Raghunathan&Etal:2008INFOCOM,Jaramillo&Srikant&Ying:2011SAC}   are also significantly different from ours.  The results in \cite{Raghunathan&Etal:2008INFOCOM,Jaramillo&Srikant&Ying:2011SAC} do not apply directly here.}

The scheduling of charging multiple EVs has received much
recent attention. In \cite{Yu&Chen&Tong:16CSEE}, the authors proposed an intelligent energy management system for the large-scale public charing stations taking into account of EV admissions, scheduling and renewable energy.
Applications of deterministic deadline scheduling models are applied
 in \cite{Chen&He&Tong:11Allerton,Chen&Ji&Tong:12PESGM,Chen&Tong:12SGC,Subramanian&Etal:12ACC}
 to study the scheduling of EV charging.
With an objective of minimizing the
load variance, a few recent papers propose several
approaches for EV charging scheduling, based on game theoretic analysis \cite{MCH10,KH13} and decentralized
algorithms \cite{Gan&Etal:2013TPS}. In addition, the authors of \cite{Benetti&Etal:TSG2015,Deilami&Etal:TSG2011,Luo&Chan:IET2014} developed control algorithm to minimize the power losses and improve the voltage profile.
 Distributed pricing strategy and algorithm are proposed in \cite{Wu&Etal:2012TSG,Donadee&Ilic:14TSG} to incentivize EVs to participate in frequency regulation. In \cite{Vagropoulos&Bakirtzis:13TPS} and \cite{Sortomme&El-Sharkawi:12TSG}, two-settlement centralized control algorithms are proposed: charging trajectories of EVs are optimized day ahead and adjustment is carried out in real-time. The authors of \cite{Juul&Etal:2015PESGM} further investigate the real-time adjustment balancing of predetermined charging trajectories according to regulation signals.

Closely related to this work,
the authors of \cite{Xu&Pan:2012CDC} construct
a dynamic framework on EV charging that explicitly
takes into account the stochasticity in both EV arrival and charging cost.
Through a dynamic programming approach,
they establish the Less Laxity and Longer Processing time (LLLP) principle:
priority should be given to vehicles with less laxity and longer processing time.
The LLLP principle is shown to be able to improve any charging policy on a sample-path basis \cite{Xu&Pan:2012CDC},
and will be used in this paper to improve the Whittle's index policy.


\section{Problem Formulation} \label{sec:II}
We now formulate the EV charging problem as a stochastic deadline scheduling problem
 subject to processing capacity  constraints.
  In Section \ref{sec:IIA}, we formulate a constrained Markov decision process (MDP).
 In Section \ref{sec:IIB}, we provide an upper bound on the total discounted reward, which  is useful for benchmark comparisons.


\begin{figure}
\centering
\psfrag{scheduling}{{\small scheduling}}
\psfrag{charging}{{\small Charging}}
\psfrag{algorithm}{{\small algorithm}}
\psfrag{cost}{{\small cost}}
\psfrag{t}{{\small $t$}}
  \includegraphics[width=.5\textwidth]{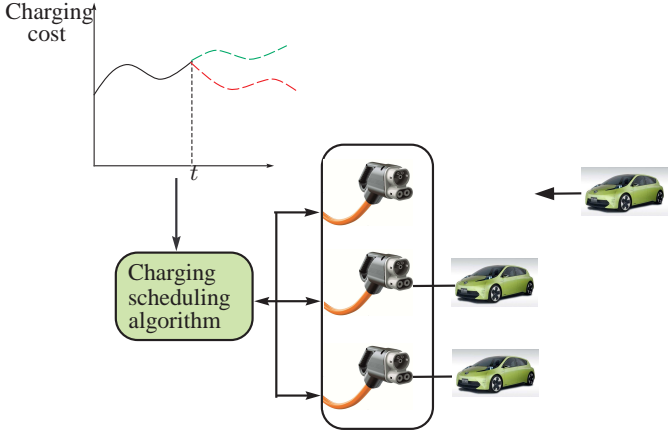}
  \caption{Architecture of a charging station}
  \label{fig:iEMS}
\end{figure}

 \subsection{An MDP Formulation of Stochastic Deadline Scheduling}\label{sec:IIA}
Fig.~\ref{fig:iEMS} shows a schematic of an energy management system at an EV charging facility.  We assume that the facility has $N$ parking spots, each with a charger  that can be activated or deactivated by the scheduler. Each charger can only be connected to one vehicle.

EVs arrive at the facility independently.  If at least one charger is available, a newly arrived EV will
park at a spot and attach to its charger.
 The EV owner communicates the charging demand $B_i$ (measured in charging time), and the deadline for completion $d_i$ to the scheduler. The scheduler receives the information and updates the state of chargers in the system. 

 %

\begin{figure}
\centering
\psfrag{time}{{\scriptsize time}}
\psfrag{t}{{\scriptsize$t$}}
\psfrag{car1}{{\scriptsize Charger $i$}}
\psfrag{rn}{{\scriptsize $r_i$}}
\psfrag{Tn}{{ \scriptsize $T_i[t]$}}
\psfrag{Bn}{{\scriptsize$B_i[t]$}}
\psfrag{xn}{{\scriptsize$L_i[t]$}}
\psfrag{dn}{{\scriptsize $d_i$}}
  \includegraphics[width=.5\textwidth]{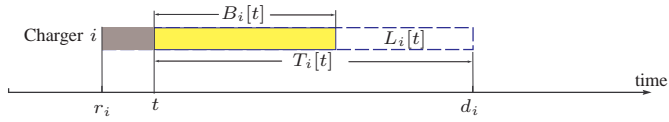}
  \caption{An illustration for the charger's state.  $r_i$ is the arrival time of an EV at charger $i$, $d_i$ the deadline for completion, $B_i[t]$ the amount of charging to be completed by $d_i$, $T_i[t]$ the lead time to deadline.}
  \label{fig:carExample}
\end{figure}

We now present elements of the discounted infinite-horizon discrete-time MDP. At the beginning of each time slot, the real time marginal charging cost is released and new EVs may arrive at the facility. The scheduler collects the states of the EVs in the facility and the charging cost, and makes a decision on which chargers to activate or deactivate in the current slot.

The assumptions in the paper are summarized as follows; they are approximations of practical operating conditions and are made for tractable analytical developments.
\bitem
\item[A1.] Each charger can be connected to only one EV, and it is removed from the EV at the deadline $d_i$.
Each EV is charged at  a fixed rate normalized to  $1$.
\item[A2.] The EV arrivals are independent.
\item[A3.] The price of charging is in proportion to the charging demand, normalized to $1$ dollar/hour.
\item[A4.] The marginal charging cost $c[t]$ is an exogenous finite state Markov chain \cite{Kwon&Xu&Gautam:TSG15} whose evolution is independent of  the state evolution and actions of charging. 
    \item[A5.] The charging of EVs is preemptive without cost.
\item[A6.] The penalty for incomplete charging  is a convex function of the incomplete amount at the deadline.
\eitem

\subsubsection{State space}
The state of the charging system consists of the state of individual chargers, charging cost, and "period index". The state of charger ${i\in\{1,\cdots,N\}}$ is defined by ${S_i[t]\defeq (T_i[t],B_i[t])}$; as illustrated in Fig.~\ref{fig:carExample}, ${T_i[t]\triangleq d_i-t}$ is the lead time and $B_i[t]$ is the remaining charging demand measured in charging time. If there is no EV attached to charger $i$, we set $S_i[t]=(0,0)$.

The system charging cost $c[t]$ is the cost of electricity. It is the electricity price from the wholesale market biased by the local renewable generation. Both of the wholesale price and the renewable generation are random and we assume the charging cost follows a Markov chain that is known to the scheduler.

We assume that the distribution of EV arrivals is time varying but periodic, i.e., the arrival distribution at the same period of each day is the same.  Such a model allows us to incorporate a ``typical day'' travel pattern for arrival statistics and convert a non-stationary arrival  to a cyclostationary one by introducing a ``period index''. Specially, each day is equally divided into $N_\tau$ periods (for example, 24 hours) and the period state $\tau[t]=(t \mbox{ mod } N_\tau)$ is the period index which forms a deterministic periodic Markov chain. The arrival rate and the probability mass function (PMF) of the initial state of EVs within the same period (for example, 9 AM-10 AM) across different days are assumed the same. The arrival rate and the initial state PMF within period $\tau$, which are known to the scheduler, are denoted by $\rho_i^\tau$ and $Q_i^\tau(T,B)$. Note that the Markov chain of charging costs needs not to be periodic (with length $N_\tau$).

Thus the state space of the charging system  is defined by
\centerline{$S[t]=(\tau[t],c[t], S_1[t],\cdots, S_N[t]) \in {\cal S}_\tau \times {\cal S}_c \times  {\cal S}_1\times \cdots \times {\cal S}_N$.}\\[0.5em]
   Here,
   ${\cal S}_\tau$ is the period space, ${\cal S}_c$ the state space of the cost, and $\mathcal{S}_i$ the state space of individual charger $i$, for $i=1,\ldots,N$.  


We note that the constructed MDP is stationary because the time dependency is incorporated by including in the system state a periodic Markov chain that describes time evolution.

\subsubsection{Action}
The action of the scheduler is defined by $\textbf{a}[t]=(a_1[t], \cdots, a_N[t]) \in \{0,1\}^N$ where $a_i[t]=1$ means that the charger is activated (active) whereas $a_i[t]=0$ means that the charger is deactivated (passive).

 \subsubsection{State evolution}
 We assume that the charging cost $c[t]$ evolves as an exogenous finite state Markov chain with transition probability matrix ${P=[P_{i,j}}]$.
The evolution of the charging cost is
 independent of the actions taken by the scheduler.

 The period state forms a deterministic periodic Markov chain. The evolution is stated as:
$ (\tau[t+1] \mid \tau[t]=\tau)=\{\tau+1\}$
where ${\{\tau+1\}=(\tau+1 \mbox{ mod }N_\tau)}$.

 Given the period index $\tau[t]=\tau$, the evolution of chargers' states depends on the scheduling action $\textbf{a}[t]=\{a_i[t]\}_{i=1}^N$, i.e.,
 \[
 (T_i[t+1],B_i[t+1]) = (T_i[t]-1, B_i[t]-a_i[t]).
 \]


 EVs leave the chargers at their deadlines.
  New EVs arrive at the charging facility following a geometric distribution with probability $\rho_i^\tau$.
   The probability mass function (PMF) $Q_i^\tau(\cdot,\cdot)$ governs the initial states of newly arrived EVs.
   The probability distribution of EV arrivals follows the periodic cycle with length $N_\tau$,
   and depends only on period index $\tau[t]$ (not $t$).
 Formally, the state evolution of charger $i$  with state $S_i[t]$ under action $a_i[t]=1$ is given by
\begin{equation}
\begin{array}{l}
\mathrel{\phantom{=}}\Big( S_i[t+1]\mid a_i[t]=1, \tau[t]=\tau\Big)\\[3pt]
=\left\{
   \begin{array}{lll}
   (T_i[t]-1,B_i[t]-1)& \mbox{w.p. } 1, \quad \mbox{if } \, B_i[t]>0, T_i[t]>1,\\[2pt]
   (T_i[t]-1,B_i[t])   & \mbox{w.p. } 1, \quad \mbox{if } \, B_i[t]=0, T_i[t]>1, \\[2pt]
   (0,0) & \mbox{w.p. } (1-\rho_i^\tau),   \quad \mbox{if } \, T_i[t]\le 1,\\[2pt]
   (1,1) & \mbox{w.p. } \rho_i^\tau Q_i^\tau(1,1), \quad \mbox{if } \, T_i[t]\le 1,\\[2pt]
   \cdots\\[2pt]
   (\bar{T},\bar{B}) & \mbox{w.p. } \rho_i^\tau Q_i^\tau(\bar{T},\bar{B}), \quad \mbox{if } \, T_i[t]\le 1, \\
   \end{array}
   \right.
   \end{array}
   \label{eqn:stateTransitionActive}
\end{equation}
where $\bar{T}$ and $\bar{B}$ is the maximum deadline and charging demand.




\subsubsection{Reward}
At time $t$, the reward received from charger $i$ under action $a_i[t]$ is given by
\beq
\begin{array}{l}
   \mathrel{\phantom{=}}R_{a_i[t]}(S_i[t],c[t])\\
   =\left\{
   \begin{array}{l}
    (1-c[t])a_i[t], \qquad \qquad \mbox{if }  B_i[t]>0, \; T_i[t]>1,\\[2pt]
   (1-c[t])a_i[t]-  F(B_i[t]-a_i[t]), \\
    \qquad \qquad\qquad \qquad \qquad \mbox{if } B_i[t]>0, \; T_i[t]=1,\\[2pt]
   0,~~~\qquad\qquad \qquad \qquad\mbox{otherwise},\\
   \end{array}
   \right.
   \label{eqn:rewardFunction}
   \end{array}
  \end{equation}
  where $F(B)$ is
  an increasing and convex penalty function with $F(0)=0$.
  Note that the scheduler obtains one unit of reward if the EV is charged for one period.
  At the EV's deadline, i.e., when $T_i[t]=1$, the scheduler pays the compensation for unfulfilled charging
  request, $F(B_i[t]-a)$.


 Given the initial system state ${S[0]=s}$ and a policy $\pi$ that maps
 each system state $S[t]$ to an action vector $\textbf{a}[t]$,
 the expected discounted system reward is defined by
 \beq \label{eq:V}
G_\pi(s) \, \defeq \,\mathbb{E}_\pi \left(\sum_{t=0}^\infty \sum_{i=1}^N \beta^t R_{a_i[t]}(S_i[t],c[t]) \mid  S[0]=s \right),
 \eeq
 where $\mathbb{E}_\pi$ is the conditional expectation over the randomness in costs and EVs arrival under a given scheduling policy $\pi$ and $0<\beta<1$ is the discount factor. The analysis can be extended to the average case \cite{Dutta:91JOE}.

 \subsubsection{Constrained MDP and optimal  policy}
We impose a constraint on the number of  simultaneously activated chargers, {\it i.e.}, ${\sum_i^N a_i[t] \le M}$ for all $t$.  In practice, such a constraint limits the peak power consumption of the  charging facility, due to feeder and transmission line capacity constraints.

The EV charging scheduling problem can be formulated as a constrained MDP.  The maximum expected reward is  given by
\beq
\label{eqn:originPro}
G(s)=\sup_{\{\pi: \sum_i^N a^{\pi}_i[t] \le M, \, \forall t\}}  G_\pi(s),
\eeq
where $a^{\pi}_i[t]$ is the action generated by policy $\pi$.  A policy $\pi^*$ is optimal if $G_{\pi^*}(s)=G(s)$.
Without loss of optimality, we will restrict our attention to  stationary policies \cite{Altman:2004book}.

\subsection{A Performance Upper Bound}\label{sec:IIB}
In  (\ref{eqn:originPro}), the power limit must be satisfied for all $t$. By relaxing this constraint
and requiring that
the average power usage does not exceed $M$, we obtain a performance upper bound for (\ref{eqn:originPro}).  In particular, a relaxed problem can be stated as
\begin{equation}
  \begin{array}{ll}
\sup_{\pi} & \mathbb{E}_\pi\left\{\sum_{t=0}^{\infty}\sum_{i=1}^{N}\beta^{t}R_{a_i[t]}(S_i[t],c[t])\mid S[0]\right\}\\[4pt]
  \mbox{subject to}& (1-\beta)\mathbb{E}\sum_{t=0}^{\infty}\sum_{i=1}^N\beta^ta_i[t]\le M.
  \end{array}
  \label{eqn:relaxedPro}
  \end{equation}
  \normalsize
Problem (\ref{eqn:relaxedPro}) is not a practical formulation for the large scale EV charging since the power usage could be far more than $M$ at certain time.

  Since the charging cost is the same for all chargers, the relaxed problem (\ref{eqn:relaxedPro}) is equivalent to
  the following problem (on the scheduling of a single charger $i$).
 \begin{equation}
  \begin{array}{ll}
  \sup_{\pi} & N\mathbb{E}_\pi\left\{\sum_{t=0}^{\infty}\beta^{t}R_{a_i[t]}(S_i[t],c[t])\mid S_i[0],c[t]\right\}\\[4pt]
  \mbox{subjec to}& (1-\beta)\mathbb{E}\sum_{t=0}^{\infty}\beta^ta_i[t]\le M/N.
  \end{array}
  \label{eqn:individualPro}
  \end{equation}
  \normalsize

  Problem (\ref{eqn:individualPro}) seeks to maximize the discounted reward from a single charger $i$ with no more than $M/N$  active action (per time period) on average.
  The optimal solution  and the optimal objective of (\ref{eqn:individualPro}) are the same as those of (\ref{eqn:relaxedPro}).
   The optimal objective of (\ref{eqn:individualPro}) can be used as a performance upper bound for the original scheduling problem in (\ref{eqn:originPro}).

The constrained MDP problem in (\ref{eqn:individualPro}) has a much smaller dimensionality and can be easily solved by linear programming (cf. Chap. 3 of \cite{Altman:2004book}
for a survey).

\section{Whittle's Index and Conservation law}\label{sec:III}
Since the MDP formulation does not result in a scalable optimal scheduling policy,
we seek to obtain an {\it index policy} \cite{Gittins:79JRSS} that can provide a scalable solution.   An index policy schedules the
  charging of EVs based on the ranked order of indices associated with the states of chargers.
  Specifically, the index of charger $i$ is a mapping from its extended  state  $\tilde{S}_i[t]\triangleq(S_i[t],c[t],\tau[t])$  to an index value. The index value of each state is independent from the states of other chargers and can be computed off-line which makes the index policy scalable.

\subsection{A Restless MAB Problem}
We now formulate problem (\ref{eqn:originPro}) as  a restless  multi-armed bandit (MAB) problem.
The restlessness is due to the fact that the state of a charger, in particular, the lead time evolves even if the charger is not activated.

A complication of casting  (\ref{eqn:originPro}) as a restless MAB problem comes from the inequality constraint on the maximum number of simultaneous activated chargers.
This complication can be circumvented by introducing  $M$ dummy chargers and requiring that exactly $M$ chargers must be activated in each period.  Specifically, each dummy charger always accrues zero reward, and the state of dummy chargers stays at ${S_i=(0,0)}$. We let $\{1,\cdots,N\}$ be the set of regular chargers and $\{N+1,\cdots,N+M\}$ be the set of
dummy chargers.

\subsubsection{Arms}
 The formulated restless  multi-armed bandit (MAB) problem has $N+M$ arms:
 each arm represents a (regular or dummy) charger.
  We define the extended state of each charger as ${\tilde{S}_i[t]\triangleq(S_i[t],c[t],\tau[t])}$, and denote the extended state space as ${\tilde{\mathcal{S}}_i\triangleq\mathcal{S}_i\times \mathcal{S}_c\times \mathcal{S}_\tau}$. The actions and the reward functions remain unchanged.


\subsubsection{MAB formulation}
By including dummy chargers, the MDP in (\ref{eqn:originPro}) is equivalent to a restless MAB problem where exactly $M$ out of ${N+M}$ chargers (arms) are active in each period. The restless MAB problem is formulated as follows
\begin{equation}
\label{eqn:MAB}
  \begin{array}{ll}
\sup_{\pi} & \mathbb{E}_\pi\left\{\sum_{t=0}^{\infty}\sum_{i=1}^{N+M}\beta^{t}R_{a_i[t]}(\tilde{S}_i[t])\mid \tilde{S}_i[0]\right\}\\[5pt]
\mbox{s.t.}& \sum_{i=1}^{N+M}a_i[t]= M, \quad {\forall } \, t.
  \end{array}
  \end{equation}
In (\ref{eqn:MAB}), the arms are coupled by the charging cost and period index, and are not independent.

\subsection{The Whittle's Index}

We now examine the Whittle's index policy for the restless MAB problem defined in  (\ref{eqn:MAB}).  To this end, we first introduce the Whittle's index and establish  the indexability of the restless MAB problem in Theorem \ref{thm:index}.

 We consider the following single charger reward maximizing problem without constraint:  given the initial state $\tilde{S}_i[0]$, policy $\pi$ activates and deactivates a single charger to maximize the reward without any power limit:
\begin{equation}V_i(\tilde{s})\triangleq\sup_{\pi}\mathbb{E}_\pi\left\{\sum_{t=0}^{\infty}\beta^{t}R_{a_i[t]}(\tilde{S}_i[t])\mid \tilde{S}_i[0]=\tilde{s} \right\},
  \label{eqn:noConstraint}
  \end{equation}
  where $V_i$ is the value function of charger  $i$. Note that
  the value function defined above is
  different from the value function defined in (\ref{eqn:originPro}) (for the constrained MDP).

Let $\mathcal{L}_{a}$ be the Markov transition operator on the extended state $\tilde{S}_i$ and an arbitrary function $f(\tilde{S}_i)$ defined as
\[
(\mathcal{L}_{a}f)(\tilde{s})
\triangleq\mathbb{E}\{f(\tilde{S}_i[t+1])\mid \tilde{S}_i[t]=\tilde{s},a_i[t]=a\}.
\]
The maximum discounted reward of problem (\ref{eqn:noConstraint}) is determined by the Bellman equation
\[
V_i(\tilde{s})=\max\{R_0(\tilde{s})+\beta(\mathcal{L}_{0}V_i)(\tilde{s}), R_1(\tilde{s})+\beta(\mathcal{L}_{1}V_i)(\tilde{s})\}.
\]
The Whittle's index is defined by introducing a  {\it $\nu$-subsidy problem},
which is a modified version of the single charger problem defined in (\ref{eqn:noConstraint}).
In the $\nu$-subsidy problem, whenever the passive action is taken, the scheduler receives an extra reward $\nu$  \cite{Whittle:1988JAP}.
The Bellman equation for the single charger $\nu$-subsidy problem is
given by
\beq\label{eqn:BellmanNuSubsidy}
V_i^\nu(\tilde{s})=\max\{R_0(\tilde{s})+\nu+\beta(\mathcal{L}_{0}V_i^\nu)(\tilde{s}), R_1(\tilde{s})+\beta(\mathcal{L}_{1}V_i^\nu)(\tilde{s})\},
\eeq
where $V_i^\nu$ is the value function for the $\nu$-subsidy problem.

Let $\tilde{\mathcal{S}}_i(\nu)$ denote the set of charger states under which it is optimal to take the passive action on charger $i$ in the $\nu$-subsidy problem. Thus any state $\tilde{s}\in\tilde{\mathcal{S}}_i(\nu)$ makes the first term in (\ref{eqn:BellmanNuSubsidy}) larger or equal to the second term.  We are now ready to define the indexability of an MAB problem.

\begin{Def}[Indexability \cite{Whittle:1988JAP}]
Charger (arm) $i$ is indexable if the set $\tilde{\mathcal{S}}_i(\nu)$ increases monotonically from $\emptyset$ to $\tilde{\mathcal{S}}_i$ as $\nu$ increases from $-\infty$ to $+\infty$. The MAB problem is indexable if all the chargers (arms) are indexable.
\end{Def}

Given the definition of indexability, the Whittle's index is defined as follows.
\begin{Def}[Whittle's index \cite{Whittle:1988JAP}]\label{def:W}
If charger (arm) $i$ is indexable, its Whittle's index $\nu_i(\tilde{s})$ of the extended state $\tilde{s}$ is the infimum subsidy $\nu$ under which the passive action is optimal at state $\tilde{s}$,
{\it i.e.},
\[
\begin{array}{l}
\nu_i(\tilde{s})\triangleq\inf_\nu\{\nu:R_0(\tilde{s})+\nu+\beta(\mathcal{L}_{0}V_i^\nu)(\tilde{s})\\[3pt]
\mathrel{\phantom{\nu_i(\tilde{s})\triangleq\inf_\nu\{\nu:}}\ge R_1(\tilde{s})+\beta(\mathcal{L}_{1}V_i^\nu)(\tilde{s})\}.
\end{array}
\]
\end{Def}
If the charger is indexable, any $\nu<\nu_s(\tilde{s})$ will make the first term strictly smaller than the second term in (\ref{eqn:BellmanNuSubsidy}) and it is optimal to activate the charger. Any $\nu\ge\nu_s(\tilde{s})$ will make the first term greater or equal to the second term and the optimal action is to deactivate the charger.

Given the definition of Whittle's index, the Whittle's index policy is stated as follows.
\begin{Def}[Whittle's index policy]\label{def:WhittleIndexPolicy}
For a multi-charger (arm) problem defined in (\ref{eqn:MAB}), the Whittle's index policy sorts all chargers by the Whittle's index value in a descend order and activates the first $M$ chargers.
\end{Def}


\subsection{Indexability and Close-form Expression of Whittle's Index}\label{sec:indexability}

In this subsection, we show that the MAB problem is indexable and the Whittle's index policy is optimal when the power limit is loose (${M=N}$). We also give the closed-form expression for the Whittle's index when the charging cost is constant.

\begin{thm}[Indexability, optimality and closed-form indexes]\label{thm:index}${}$\\
$1)$ Each charger is indexable and the MAB problem (\ref{eqn:MAB}) is indexable.\\
$2)$  When ${M=N}$, the Whittle's index policy is optimal for the multi-armed bandit problem defined in (\ref{eqn:MAB}).\\
$3)$ If $c[t]=c_0$ for all $t$, the Whittle's index of a regular charger $i\in\{1,\cdots,N\}$ is given by
\begin{equation}\label{eqn:closedForm}
\hspace{-1em}
\begin{array}{l}
\mathrel{\phantom{=}}\nu_i(T,B,c_0,\tau)\\
=\left\{
\begin{array}{ll}
0 &\mbox{if~}\, B=0,\\[3pt]
1-c_0 &\mbox{if~}\, 1\le B\le T-1,\\[3pt]
1-c_0+\\
\beta^{T-1}[F(B-T+1)-F(B-T)]&\mbox{if~}\, T\le B.
\end{array}
\right.
\end{array}
\end{equation}
The Whittle's index of a dummy charger is zero.
\[
\nu_i(0,0,c_0,\tau)=0, \quad i\in\{N+1,\cdots,N+M\}.
\]

\end{thm}

\begin{proof}
 An elementary proof of indexability can be found in Appendix \ref{proof:indexability}. The proof of optimality of Whittle's index policy with $M=N$ can be found in Appendix \ref{proof:optimalityOfWhittleIndexPolicy}. The proof of closed-form of Whittle's index with constant charging cost can be found in Appendix \ref{proof:closed-form}.

\end{proof}

In (\ref{eqn:closedForm}), when it is feasible to fulfill  EV $i$'s charging request ({\it i.e.} its lead time is no less than its remaining processing time),
EV $i$'s Whittle's index is simply the (per-unit) charging profit $1-c_0$. When non-completion penalty is inevitable, the index takes into account both
the charging profit and the non-completion penalty. We note that the Whittle's index gives higher priority to EVs with less laxity.
Here, the laxity of charger $i$ is defined as $L_i[t]\triangleq T_i[t]-B_i[t]$ (cf. Fig.~\ref{fig:carExample}).

When the power limit is loose (${M=N}$), the MAB problem breaks into $N$ independent single arm problems and Whittle's index policy is optimal. It balances the charging cost and the penalty of non-finished demand by deactivating (regular) chargers when the charging cost is high. Simple index policies such as the earliest deadline first (EDF) and least laxity first (LLF) policies do not take charging cost into account and may lead to significant performance loss.
However, when $M<N$, we note that the Whittle's index policy does not distinguish EVs with positive laxity, and is therefore suboptimal. In the next section we will introduce an enhanced heuristic policy based on the Whittle's index.

\section{Whittle's Index Policy with LLLP Interchange}\label{sec:IV}
For the objective of time average (${\beta=1}$) profit maximization,
the Whittle's index policy is shown (under some conditions on the evolution of arm states) to be  asymptotically optimal, as the number of arms increases to infinity \cite{Mora:2001AAP}.
For the discounted profit maximization setting considered in this paper,
the asymptotic optimality of Whittle's index policy is not clear.
For small systems with finitely many arms, there are counter-examples where an optimal index policy does not exist (and therefore the Whittle's index policy cannot be optimal).

In this section, we will apply the Less Laxity and Longer remaining Processing time (LLLP) principle (originally proposed in \cite{Xu&Pan:2012CDC}) to improve the Whittle's index policy.

The LLLP principle is a priority rule for the scheduling
of charging multiple EVs, which is defined as follows.
\begin{Def}[The LLLP Principle] \label{def:LLLP}
Consider chargers (arms) $i$ and $j$ at time $t$. We say $j$ dominates $i$ ($j\succeq i$),
if $j$ has less laxity and longer remaining processing time, \ie  $L_j[t]\le L_i[t]$ and $B_j[t] \ge B_i[t]$, with at least one of the inequalities strictly holds.
\end{Def}

LLLP defines a partial order over the EVs' states such that the EV with less laxity and longer remaining charging demand should be given priority.
In \cite{Xu&Pan:2012CDC}, the authors applied an interchange argument to show that LLLP could improve the performance of any given policy along every sample path,
 and further, there exists an optimal stationary policy that follows the LLLP principle under mild conditions.

To apply the LLLP principle, note that the Whittle's index policy for the multi-armed bandit problem is a stationary policy:
at each time it orders (the states of) the $M+N$ arms, and activates
the first $M$ arms. The proposed heuristic policy re-order
every pair of arms that violates the LLLP principle (cf. Algorithm~\ref{alg:WhittleLLLP}).
As such, the proposed heuristic policy always gives priority to EVs with less laxity and longer remaining processing time.



\begin{algorithm}
\caption{Whittle index with LLLP interchange}
\label{alg:WhittleLLLP}
\begin{algorithmic}
\STATE $1$. Calculate the Whittle's indexes of all chargers and sort them in a descend order.
\STATE $2$. Suppose the order of chargers is $i_1,i_2,\cdots,i_{M+N}$.\\
for $j=i_1:i_M$\\
~~~for $k=i_{M+1}:i_{M+N}$\\
~~~~~~~if $k\succeq j$ in the sense of Definition \ref{def:LLLP}\\
~~~~~~~~~~exchange the orders of $j$ and $k$\\
~~~~~~~end\\
~~~end\\
end
\STATE $3$. Activate the $M$ chargers with highest priority.
\end{algorithmic}
\end{algorithm}

\section{Numerical Results}\label{sec:V}
\subsection{Benchmark Policies}
In this section, results of numerical experiments are presented to compare  the performance of proposed index policy with other simple heuristic (index) policies,
 {\it i.e.}, EDF (earliest deadline first) \cite{Liu&Layland:73ACM}, LLF (least laxity first) \cite{Dertouzos:74IFIPC}, valley filling \cite{Gan&Etal:2013TPS}, and the original Whittle's index (without LLLP interchange) \cite{Whittle:1988JAP}.

 If feasible, EDF charges $M$ EVs with the earliest deadlines, and  LLF charges $M$ EVs with the least laxity.
Both policies will fully utilize the capacity and charge $M$ EVs as long as there are at least $M$ unfinished EVs in the system.
The Whittle's index policy, on the other hand,  ranks all chargers by the Whittle's index and activates the first $M$ arms, and may
put some (regular) chargers idle (deactivated) when the charging cost is high.

The centralized valley filling algorithm  is proposed in \cite{Gan&Etal:2013TPS} as an optimal offline scheduling policy for a setting with continuous charging rate, deterministic charging cost, and no newly arrival EVs. In our setting with newly incoming EVs, the algorithm is repeatedly executed  based on
  the most updated information on the arrived EVs. At the beginning of each period, the algorithm
   schedules EV charging based on the information of all EVs that have arrived, assuming there is no new incoming EVs in the future. The charging schedule is carried out only for current time slot; at the beginning of the next period, the scheduling plan is recomputed.
   The (real-time adjusted version of) valley filling algorithm is used in our simulation as a benchmark.\footnote{The original valley filling algorithm proposed in \cite{Gan&Etal:2013TPS} has a hard constraint that all EVs' charging request must be fulfilled. We note that it may not be feasible to fulfill all EVs' charging requests
under the charging capacity $M$. When implementing the valley filling algorithm in our setting, we modify the objective function of
 the
 valley filling algorithm so as to maximize the total reward subject to the maximum power limit; in its objective function, the charging reward is calculated based on expected future charging cost and the non-completion penalty function $F(\cdot)$  (introduced in (\ref{eqn:rewardFunction})).}


\subsection{Numerical Results}
We first considered a special case of problem (\ref{eqn:MAB})
  with a constant charging cost.
    Since the charging cost was time-invariant,
    it was optimal to fully utilize the charging capacity to charge $M$ unfinished EVs.

\begin{figure}[h]
\centering
\psfrag{Num of arm}{{\small $M/N$}}
\psfrag{Reward per arm}{{\small Total reward/\# of regular chargers (\$)}}
\psfrag{UB}{{\tiny Upper Bound}}
\psfrag{EDF}{{\tiny EDF}}
\psfrag{LLF}{{\tiny LLF}}
\psfrag{Valley filling}{{\tiny Valley filling}}
\psfrag{Whittls index}{{\tiny Whittle's index}}
\psfrag{LLLLLLLLLLLLP-Index}{{\tiny Whittle's index w. LLLP}}
  \includegraphics[width=.45\textwidth]{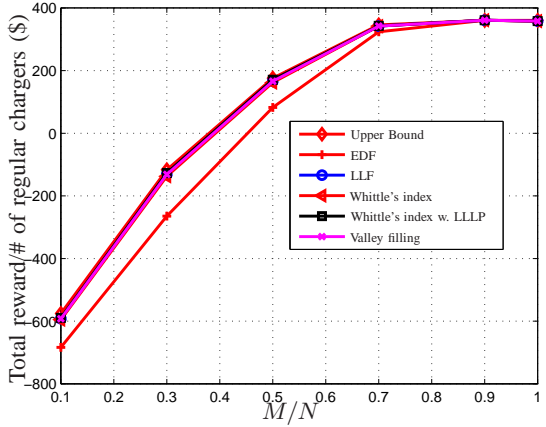}
  \caption{Performance comparison: constant charging cost ${c[t]=0.5}$, ${\rho_i^\tau=0.7}$, ${\bar{T}=12}$, ${\bar{B}=9}$, ${\beta=0.999}$, ${F(B)=0.2B^2}$, ${N=10}$.}
  \label{fig:constantCostVaryRatio}
\end{figure}

\begin{figure}[h]
\centering
\psfrag{Num of arm}{{\small \# of regular chargers}}
\psfrag{Reward per arm}{{\small Total reward/\# of regular chargers (\$)}}
\psfrag{UB}{{\tiny Upper Bound}}
\psfrag{EDF}{{\tiny EDF}}
\psfrag{LLF}{{\tiny LLF}}
\psfrag{Valley filling}{{\tiny Valley filling}}
\psfrag{Whittl's index}{{\tiny Whittle's index}}
\psfrag{LLLLLLLLLLLLP-Index}{{\tiny Whittle's index w. LLLP}}
  \includegraphics[width=.45\textwidth]{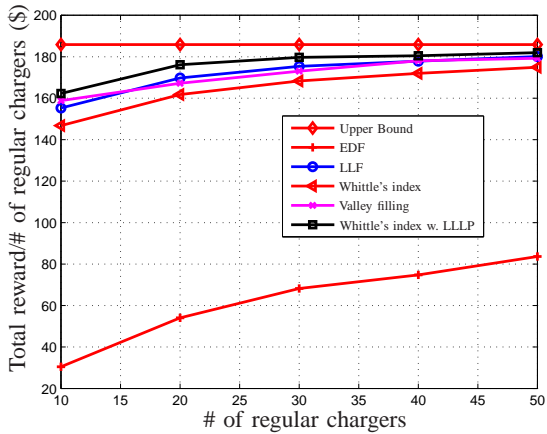}
  \caption{Performance comparison: constant charging cost ${c[t]=0.5}$, ${\rho_i^\tau=0.7}$, ${\bar{T}=12}$, ${\bar{B}=9}$, ${\beta=0.999}$, ${F(B)=0.2B^2}$, ${M/N=0.5}$.}
  \label{fig:constantCost}
\end{figure}

In Fig.~\ref{fig:constantCostVaryRatio}, we fixed the traffic of EVs and the total number of regular chargers and varied the power limit $M$. All policies besides the EDF scheduling performed well and close to the upper bound of the performance. When $M/N=1$, all EVs can be fully charged. Thus all policies achieved optimality. In Fig. \ref{fig:constantCost}, we zoomed in the case when $M/N=0.5$ and varied the total number of regular charger $N$. We observed that the Whittle's index policy with LLLP interchange and LLF achieved similar performance, since both policies roughly followed the least laxity first principle. The performance of these two policies was close to the performance upper bound.
The EDF policy performed poorly because it did not take the remaining charging demand into account.
The gap between the Whittle's index policy and the Whittle's index policy with LLLP interchange came from the reordering of EVs with positive laxity (cf. the discussion following
Theorem \ref{thm:index}).

For the dynamic charging cost case, we used the real-time pricing signal from the California Independent System Operator (CAISO) and trained a Markovian model that describes the
marginal charging costs (cf. Sections III and V of \cite{Kwon&Xu&Gautam:TSG15}). Each period of the
constructed Markov chain (on charging cost) lasted for 1 hour, and
each periodic cycle lasted for one day with ${N_\tau=24}$.
For each period, the real-time price was quantized into discrete price states, and the transition probability (of the Markov chain) was simply the frequency the price changed from one state to another.

\begin{figure}[h]
\centering
\psfrag{Num of arm}{{\small $M/N$}}
\psfrag{Reward per arm}{{\small Total reward/\# of regular chargers (\$)}}
\psfrag{UB}{{\tiny Upper Bound}}
\psfrag{EDF}{{\tiny EDF}}
\psfrag{LLF}{{\tiny LLF}}
\psfrag{Valley filling}{{\tiny Valley filling}}
\psfrag{Whittls index}{{\tiny Whittle's index}}
\psfrag{LLLLLLLLLLLLP-Index}{{\tiny Whittle's index w. LLLP}}
  \includegraphics[width=.45\textwidth]{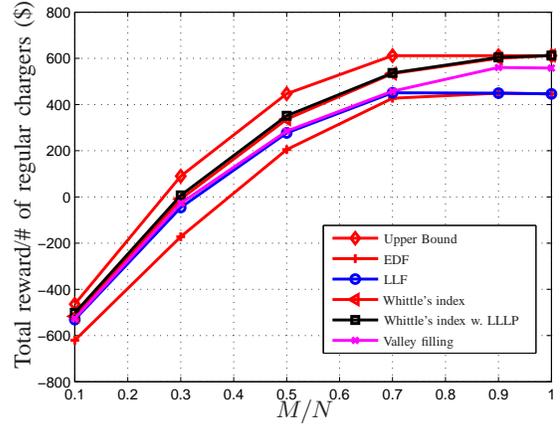}
  \caption{Performance comparison: dynamic charging cost, ${\rho_i^\tau=0.7}$, ${\bar{T}=12}$, ${\bar{B}=9}$, ${\beta=0.999}$, ${F(B)=0.2B^2}$.}
  \label{fig:varyRatio}
\end{figure}

In Fig.~\ref{fig:varyRatio}, we fixed the EV traffic and the total number of regular chargers ${N=10}$ and varied the power limit constraints. When the power limit was low and $M/N$ is small, there was no enough power to charge EVs and the penalty dominated the charging profit. In this case the performance of different policies were close since the limited resource allowed not much to do. When the power limit was adequate and ${M=N}$, all EVs could be fully charged on time. In this case, the Whittle's index policy solved the problem optimally as stated in Theorem \ref{thm:index} and achieved the upper bound. The interchange did not happen because the LLLP principle was always satisfied in this case. The valley filling algorithm did not consider the future arrivals and EDF and LLF did not consider the dynamic charging cost, thus they performed sub-optimal. When the power constraint was neither too tight ($M/N\approx0$) nor too loose ($M/N\approx1$), LLLP could reduce the number of unfinished EVs with large remaining charging demand and therefore reduced the non-completion penalties.

\begin{figure}
\centering
\psfrag{Num of arm}{{\small \# of regular chargers}}
\psfrag{Reward per arm}{{\small Total reward (\$)}}
\psfrag{Reward}{{\small Total reward}}
\psfrag{UB}{{\tiny Upper Bound}}
\psfrag{EDF}{{\tiny EDF}}
\psfrag{LLF}{{\tiny LLF}}
\psfrag{SL}{{\tiny Valley filling}}
\psfrag{Index-Idle}{{\tiny Whittle's index}}
\psfrag{Index-Idle-NoBack}{{\tiny Whittle's index w. LLLP}}
\psfrag{Whittls index}{{\tiny Whittle's index}}
\psfrag{LLLLLLLLLLLLP-Index}{{\tiny Whittle's index w. LLLP}}
\psfrag{Valley filling}{{\tiny Valley filling}}
  \includegraphics[width=.48\textwidth]{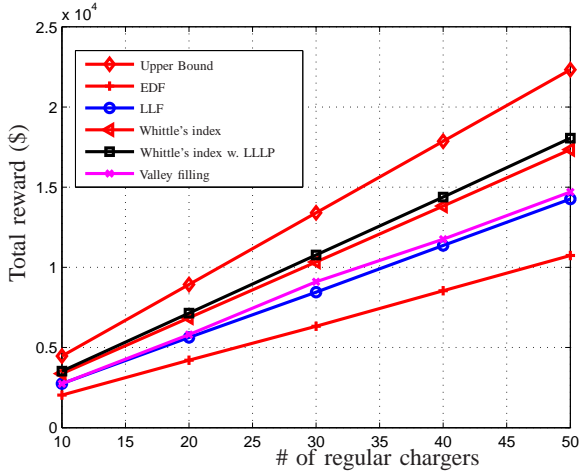}
  \caption{Performance comparison: dynamic charging cost, ${\rho_i^\tau=0.7}$, ${\bar{T}=12}$, ${\bar{B}=9}$, ${\beta=0.995}$, ${F(B)=0.2B^2}$, ${M/N=0.5}$.}
  \label{fig:dynamicCost}
\end{figure}

In Fig.~\ref{fig:dynamicCost}, we compared the performance of different policies by fixing the power limit ratio $M/N=0.5$ and varied the number of regular chargers. Both EDF and LLF sought to activate as many regular chargers as possible, up to the capacity constraint $M$.
The Whittle's index policy, on the other hand, took the advantage of the pricing fluctuation and charged more EVs at price valley and kept some chargers idle when charging cost was high. Based on the Whittle's index policy, the LLLP interchange reduced the penalty of unfinished EVs and improved the performance of Whittle's index policy.
The total reward achieved by the Whittle's index with LLLP interchange policy was more than 1.7 times of that obtained by EDF; the performance gap between the
the Whittle's index with LLLP interchange policy and the LLF policy was over $25\%$. We also noted that the LLLP principle improved Whittle's index by around $10\%$.


\section{Conclusion}\label{sec:VI}
We considered the problem of scheduling of large scale EV charging in public facilities---a problem of particular potential significance as EV penetration deepens.  In such settings, it is essential to develop highly efficient and online charging algorithms.  To this end, index policies considered in this paper are attractive for its implementation simplicity and versatility in incorporating various operation uncertainties.    
\appendix

\subsection{Proof of Indexability}\label{proof:indexability}
In this subsection, we provide an elementary proof of indexability. That is, for any charger state $\tilde{s}$, there is a critical $\nu(\tilde{s})$ such that if and only if $\nu\ge\nu(\tilde{s})$ the first term in the Bellman equation (\ref{eqn:BellmanNuSubsidy}) is larger or equal to the the second term in a single charger $\nu$-subsidy problem.

The indexability of dummy chargers is straightforward. For ${i\in\{N+1,\cdots,N+M\}}$, there is no EV arrival, and only the charging cost and period index evolve.
The Bellman equation of the $\nu$-subsidy problem is given by
\[\begin{array}{l}
V_i^\nu(0,0,c_j,\tau)=\max\{\beta\sum_kP_{j,k}V^\nu_i(0,0,c_k,\{\tau+1\})+\nu,\\
\mathrel{\phantom{V_i^\nu(0,0,c_j,\tau)=\max\{}}\beta\sum_kP_{j,k}V^\nu_i(0,0,c_k,\{\tau+1\})\}.
\end{array}\]

If and only if ${\nu\ge0}$, the first term is larger than the second term and it is optimal to deactivate the dummy charger. Otherwise, the active action is optimal. So a dummy charger is indexable and its Whittle's index is ${\nu_i(0,0,c_j,\tau)=0}$.

For the regular chargers, we prove the indexability by induction.
\subsubsection{When $T=0$} There is no EV attached to the charger. The Bellman equation is stated as
    \[V_i^\nu(0,0,c_j,\tau)=\max\{\nu+\beta W_{j,\tau}^\nu,\beta W_{j,\tau}^\nu\}.\]
        where
            \[
            \begin{array}{l}
               \mathrel{\phantom{=}} W_{j,\tau}^\nu\\
               =(1-\rho_i^\tau)\sum_{k}P_{j,k}V_i^\nu(0,0,c_k,\{\tau+1\})+\\
                \mathrel{\phantom{=}}\rho_i^\tau\sum_{T}\sum_{B}\sum_{k}Q_i^\tau(T,B)P_{j,k}V_i^\nu(T,B,c_k,\{\tau+1\})
                        \end{array}
                            \] is the expected reward of possible arrivals. If and only if ${\nu\ge0}$, the first term is larger and the passive action is optimal. Thus ${\nu_i(0,0,c_j,\tau)=0}$.

\subsubsection{When $T=1$}\label{sec:closedForm1} There are two possible conditions.
    \bitem
    \item If $B=0$, the Bellman equation is stated as
\[V_i^\nu(1,0,c_j,\tau)=\max\{\nu+\beta W_{j,\tau}^\nu,\beta W_{j,\tau}^\nu\}.\]
Thus ${\nu_i(1,0,c_j,\tau)=0}$.
    \item If $B\ge1$, the Bellman equation is stated as
    \[\hspace{-1em}\begin{array}{l}
    V_i^\nu(1,B,c_j,\tau)=\max\{\nu-F(B)+\beta W_{j,\tau}^\nu,\\
    \mathrel{\phantom{V_i^\nu(1,B,c_j,\tau)=\max\{}}1-c_j-F(B-1)+\beta W_{j,\tau}^\nu\}.
    \end{array}\]
    \normalsize
   If and only if ${\nu\ge1-c_j+F(B)-F(B-1)}$, the passive action is optimal.
\eitem

   Thus the arm is indexable and the Whittle's index for ${T=1}$ is stated as follows.
    \[
    \begin{array}{l}
    \mathrel{\phantom{=}}\nu_i(1,B,c_j,\tau)\\
    =\left\{
    \begin{array}{ll}
    0&\mbox{if~} B=0\\
    1-c_j+F(B)-F(B-1)&\mbox{if~} B\ge1
    \end{array}
    \right.\end{array}
    \]

We summarize some properties of the Whittle's index prove them for the case ${T=1}$.
    \begin{prop}[Monotonicity of Whittle's index]\label{Prop:monotonicityOfIndex1}
    Whittle's index is non-decreasing in the charging demand $B$ when $B\ge T$.
       \[\nu_i(T,B,c_j,\tau)\le\nu_i(T,B+1,c_j,\tau), \forall B\ge T.\]
       \end{prop}
       \begin{proof}
       Since the penalty function $F(\cdot)$ is convex, Whittle's index is nondecreasing when $B\ge T$ and $T=1$.
       \end{proof}

       \begin{prop}[Difference of value functions]\label{Prop:differenceOfValueFunction1}
       Denote the difference of the value function as
     \[
     \begin{array}{l}
       \mathrel{\phantom{=}}g_h^\nu(T,B,c_j,\tau)=V_i^\nu(T,B+h,c_j,\tau)-V_i^\nu(T,B,c_j,\tau),
       \end{array}
       \]
       \normalsize
        where ${h\in\{1,\cdots,\bar{B}-B\}}$ and $\bar{B}$ is the maximum charging demand.

        Thus $g_h^\nu$ is piecewise linear in $\nu$ and has following properties.
       \bitem
       \item $g_h^\nu(T,B,c_j,\tau)$ is continuous in $\nu$.

           \item There exist $\underline{\nu}_h(T,B,c_j,\tau)$ and $\bar{\nu}_h(T,B,c_j,\tau)$ such that, ${g_h^\nu(T,B,c_j,\tau)}$ is piecewise linear and ${\partial g_h^\nu/\partial\nu\ge-h}$ when ${\nu\in[\underline{\nu}_h(T,B,c_j,\tau),\bar{\nu}_h(T,B,c_j,\tau)]}$. Otherwise, $g_h^\nu(T,B,c_j,\tau)$ is constant .
               \eitem
       \end{prop}

       \begin{proof}
       \bitem
\item When $B=0$, \[g_h^\nu(1,B,c_j,\tau)=V_i^\nu(1,h,c_j,\tau)-V_i^\nu(1,0,c_j,\tau).\]

If $\nu_i(1,h,c_j,\tau)>\nu_i(1,0,c_j,\tau)=0$,
\[\hspace{-1em}
\begin{array}{l}
\mathrel{\phantom{=}}g_h^\nu(1,B,c_j,\tau)\\
=\left\{
\begin{array}{ll}
1-c_j-F(h-1),&\mbox{if~}\nu<0;\\
1-c_j-F(h-1)-\nu,&\mbox{if~}0\le\nu<\nu_i(1,h,c_j,\tau);\\
-F(h),&\mbox{if~}\nu_i(1,h,c_j,\tau)\le\nu.
\end{array}
\right.
\end{array}\]

If ${\nu_i(1,h,c_j,\tau)\le\nu_i(1,0,c_j,\tau)=0}$,
\[\hspace{-1em}\begin{array}{l}
\mathrel{\phantom{=}}g_h^\nu(1,B,c_j,\tau,)\\
=\left\{
\begin{array}{ll}
1-c_j-F(h-1),&\mbox{if~}\nu<\nu_i(1,h,c_j,\tau);\\
\nu-F(h),&\mbox{if~}\nu_i(1,h,c_j,\tau)\le\nu<0;\\
-F(h),&\mbox{if~}\nu_i(1,h,c_j,\tau)\le\nu.
\end{array}
\right.
\end{array}\]

\item When $B\ge1$, \[{g_h^\nu(1,B,c_j,\tau)=V_i^\nu(1,B+h,c_j,\tau)-V_i^\nu(1,B,c_j,\tau)}.\]
Since ${\nu_i(1,h+B,c_j,\tau)\ge\nu_i(1,B,c_j,\tau)}$,
\[\hspace{-1em}\begin{array}{l}
\mathrel{\phantom{=}}g_h^\nu(1,B,c_j,\tau)\\
=\left\{
\begin{array}{l}
F(B-1)-F(B+h-1),\\
\mbox{if~}\nu<\nu_i(1,B,c_j,\tau);\\
1-c_j-\nu+F(B)-F(B+h-1),\\
\mbox{if~}\nu_i(1,B,c_j,\tau)\le\nu<\nu_i(1,B+h,c_j,\tau);\\
F(B)-F(B+h),\\
\mbox{if~}\nu_i(1,B+h,c_j,\tau)\le\nu.
\end{array}
\right.
\end{array}\]
\normalsize
\eitem
Thus, $g_h^\nu(1,B,c_j,\tau)$ is piecewise linear and continuous in $\nu$. The derivative satisfies Property \ref{Prop:differenceOfValueFunction1}.
\end{proof}

       \begin{prop}[Concavity of value functions]\label{Prop:concavityOfValueFunction1}
       The difference of value functions ${g_1^\nu(T,B,c_j,\tau)}$ is non-increasing in $B$ when ${B\ge T}$. Thus $V_i^\nu(T,B,c_j,\tau)$ is concave in $B$ when ${B\ge T}$.
       \end{prop}
       \begin{proof}
       Since when $B\ge1$, $\nu_i(1,B,c_j,\tau)$ is nondecreasing in $B$,
\[\begin{array}{l}
\mathrel{\phantom{=}}g_1^\nu(1,B+1,c_j,\tau)-g_1^\nu(1,B,c_j,\tau)\\
=V_i^\nu(1,B+2,c_j,\tau)+V_i^\nu(1,B,c_j,\tau)-\\
\mathrel{\phantom{=}}2V_i^\nu(1,B+1,c_j,\tau)\\
=\left\{
\begin{array}{l}
2F(B)-F(B+1)-F(B-1),\\
\mbox{if~}\nu\le\nu_i(1,B,c_j,\tau);\\
\nu-(1-c_j)+F(B)-F(B+1),\\
\mbox{if~}\nu_i(1,B,c_j,\tau)<\nu\le\nu_i(1,B+1,c_j,\tau);\\
1-c_j-\nu-F(B)+F(B+1),\\
\mbox{if~}\nu_i(1,B+1,c_j,\tau)<\nu\le\nu_i(1,B+2,c_j,\tau);\\
2F(B+1)-F(B)-F(B+2),\\
\mbox{if~}\nu_i(1,B+2,c_j,\tau)\le\nu;
\end{array}
\right.\\
\le0.
\end{array}
\]
\normalsize
The first and last cases are non-positive since the penalty function $F(B)$ is convex. The second and third cases are non-positive because of the expressions of the index $\nu_i$. Thus $V_i^\nu(1,B,c_j,\tau)$ is concave in $B$ when $B\ge1$.
       \end{proof}

\subsubsection{$T\ge2$} Assuming indexability and Proposition~\ref{Prop:monotonicityOfIndex1},~\ref{Prop:differenceOfValueFunction1},~\ref{Prop:concavityOfValueFunction1} hold for $T-1$, we show that they hold for $T$.

       %
%
If $B=0$, the Bellman equation is stated as follows.
    \[
    \hspace{-1em}\begin{array}{l}
    \mathrel{\phantom{=}}V_i^\nu(T,0,c_j,\tau)\\
    =\max\{\beta\sum_kP_{j,k}V_i^\nu(T-1,0,c_k,\{\tau+1\})+\nu,\\
    \mathrel{\phantom{=\max\{}}\beta\sum_kP_{j,k}V_i^\nu(T-1,0,c_k,\{\tau+1\})\}.
    \end{array}\]
    If and only if ${\nu\ge0}$, the first term is larger than the second term and the passive action is optimal. Thus ${\nu_i(T,0,c_j,\tau)=0}$.

If $B\ge1$, the Bellman equation is stated as follows.
    \begin{equation}\label{eqn:BLargerThanOen}
    \begin{array}{l}
    \mathrel{\phantom{=}}V_i^\nu(T,B,c_j,\tau)\\
    =\max\{\beta\sum_kP_{j,k}V_i^\nu(T-1,B,c_k,\{\tau+1\})+\nu,\\
    \mathrel{\phantom{=\max\{}}\beta\sum_kP_{j,k}V_i^\nu(T-1,B-1,c_k,\{\tau+1\})+1-c_j\}.
    \end{array}
    \end{equation}
    Denote the difference between the two actions as
    \[
    \begin{array}{l}
     f^\nu(T,B,c_j,\tau)\triangleq\nu-(1-c_j)+\\
   \mathrel{\phantom{\triangleq}}\beta\sum_kP_{j,k}g_1^\nu(T-1,B-1,c_k,\{\tau+1\}),
    \end{array}
    \]
    where \[
    \begin{array}{l}
    \mathrel{\phantom{=}}g_1^\nu(T-1,B-1,c_k,\tau)\\
    =V_i^\nu(T-1,B,c_k,\tau)-V_i^\nu(T-1,B-1,c_k,\tau).
    \end{array}\]

Since Proposition~\ref{Prop:differenceOfValueFunction1} holds for ${T-1}$, ${f^\nu(T,B,c_j,\tau)}$ is continuous and piece-wise linear in $\nu$. Denote
\[\begin{array}{l}
\underline{\nu}(T,B,c_j,\tau)\triangleq\min_k\underline{\nu}_1(T-1,B-1,c_k,\tau),\\
\bar{\nu}(T,B,c_j,\tau)\triangleq\max_k\bar{\nu}_1(T-1,B-1,c_k,\tau).
\end{array}\]
We have
\[
\begin{array}{l}
\mathrel{\phantom{=}}\partial f^\nu(T,B,c_j,\tau)/\partial\nu\\
=\left\{
\begin{array}{ll}
1,& \mbox{if~}\nu\notin[\underline{\nu}(T,B,c_j,\tau),\bar{\nu}(T,B,c_j,\tau)];\\
\ge0,& \mbox{otherwise.}
\end{array}
\right.
\end{array}
\]

So $f^\nu(T,B,c_j,\tau)$ is continuous and non-decreasing in $\nu$. When ${\nu=-\infty}$, ${f^\nu(T,B,c_j,\tau)=-\infty}$. When ${\nu=+\infty}$, ${f^\nu(T,B,c_j,\tau)=+\infty}$. Thus there is a cross point of $f$ and the x-axis. Define ${\nu(T,B,c_j,\tau)\triangleq\min_{\nu}\{f^\nu(T,B,c_j,\tau)=0\}}$. If and only if ${\nu\ge\nu(T,B,c_j,\tau)}$, the first term in (\ref{eqn:BLargerThanOen}) is larger or equal to the second term and the passive action is optimal. By definition, $\nu(T,B,c_j,\tau)$ is the Whittle's index.

The indexability of $T$ is shown. Next we will prove Proposition~\ref{Prop:monotonicityOfIndex1},~\ref{Prop:differenceOfValueFunction1},~\ref{Prop:concavityOfValueFunction1} for $T$ assuming that they are true for ${T-1}$.

 \begin{proof}[Proof of Proposition~\ref{Prop:monotonicityOfIndex1}]
     When $B\ge T$,
       \[
       \begin{array}{l}
   \mathrel{\phantom{=}}    f^\nu(T,B,c_j,\tau)-f^\nu(T,B+1,c_j,\tau)\\
       =\beta\sum_kP_{j,k}g_1^\nu(T-1,B-1,c_k,\{\tau+1\})-\\
       \mathrel{\phantom{=}}\beta\sum_kP_{j,k}g_1^\nu(T-1,B,c_k,\{\tau+1\})\\
       \ge0
       \end{array}
       \]
\normalsize
       The inequality is because that $V_i^\nu(T-1,B,c_i,\tau)$ is concave in $B$ when $B\ge T-1$.

 Since $f^\nu(T,B,c_j,\tau)\ge f^\nu(T,B+1,c_j,\tau)$ we have
           \[\nu(T,B,c_j,\tau)\le \nu(T,B+1,c_j,\tau), \forall B\ge T.\]
           \end{proof}

\begin{proof}[Proof of Proposition~\ref{Prop:differenceOfValueFunction1}]

If $B=0$, \[g_h^\nu(T,B,c_j,\tau)=V_i^\nu(T,h,c_j,\tau)-V_i^\nu(T,0,c_j,\tau).\]
\bitem
\item If $\nu_i(T,h,c_j,\tau)>\nu_i(T,0,c_j,\tau)=0$,
\[\hspace{-1em}
\begin{array}{ll}
\mathrel{\phantom{=}}g_h^\nu(T,0,c_j,\tau)\\
=\left\{
\begin{array}{l}
1-c_j+\beta\sum_kP_{j,k}g_{h-1}^\nu(T-1,0,c_k,\{\tau+1\}),\\
\mbox{if~}\nu<0;\\
1-c_j-\nu+\beta\sum_kP_{j,k}g_{h-1}^\nu(T-1,0,c_k,\{\tau+1\}),\\
\mbox{if~}0\le\nu<\nu_i(T,h,c_j,\tau);\\
\beta\sum_kP_{j,k}g_{h}^\nu(T-1,0,c_k,\{\tau+1\}),\\
\mbox{if~}\nu_i(T,h,c_j,\tau)\le\nu.
\end{array}
\right.
\end{array}
\]
\item If $\nu_i(T,h,c_j,\tau)\le\nu_i(T,0,c_j,\tau)=0$,
\[\hspace{-1em}
\begin{array}{l}
\mathrel{\phantom{=}}g_h^\nu(T,0,c_j,\tau)\\
=\left\{
\begin{array}{l}
1-c_j+\beta\sum_kP_{j,k}g_{h-1}^\nu(T-1,0,c_k,\{\tau+1\}),\\
\mbox{if~}\nu<\nu_i(T,h,c_j,\tau);\\
\nu+\beta\sum_kP_{j,k}g_{h}^\nu(T-1,0,c_k,\{\tau+1\}),\\
\mbox{if~}\nu_i(T,h,c_j,\tau)\le\nu<0;\\
\beta\sum_kP_{j,k}g_{h}^\nu(T-1,0,c_k,\{\tau+1\}),\\
\mbox{if~}0\le\nu.
\end{array}
\right.
\end{array}
\]
\eitem

If $B\ge1$, \[g_h^\nu(T,B,c_j,\tau)=V_i^\nu(T,B+h,c_j,\tau)-V_i^\nu(T,B,c_j,\tau).\]
\bitem
\item If $\nu_i(T,B+h,c_j,\tau)>\nu_i(T,B,c_j,\tau)$,
\[\hspace{-1em}
\begin{array}{l}
\mathrel{\phantom{=}}g_h^\nu(T,B,c_j,\tau)\\
=\left\{
\begin{array}{l}
\beta\sum_kP_{j,k}g_h^\nu(T-1,B-1,c_k,\{\tau+1\}),\\
\mbox{if~}\nu<\nu_i(T,B,c_j,\tau);\\
1-c_j-\nu+\\
\beta\sum_kP_{j,k}g_{h-1}^\nu(T-1,B,c_k,\{\tau+1\}),\\
\mbox{if~}\nu_i(T,B,c_j,\tau)\le\nu<\nu_i(T,B+h,c_j,\tau);\\
\beta\sum_kP_{j,k}g_{h}^\nu(T-1,B,c_k,\{\tau+1\}),\\
\mbox{if~}\nu_i(T,B+h,c_j,\tau)\le\nu.
\end{array}
\right.
\end{array}
\]
\item If $\nu_i(T,B+h,c_j,\tau)\le\nu_i(T,B,c_j,\tau)$,
\[\hspace{-1em}
\begin{array}{l}
\mathrel{\phantom{=}}g_h^\nu(T,B,c_j,\tau)\\
=\left\{
\begin{array}{l}
\beta\sum_kP_{j,k}g_h^\nu(T-1,B-1,c_k,\{\tau+1\}),\\
\mbox{if~}\nu<\nu_i(T,B+h,c_j,\tau);\\
\nu-(1-c_j)+\\
\beta\sum_kP_{j,k}g_{h+1}^\nu(T-1,B-1,c_k,\{\tau+1\}),\\
\mbox{if~}\nu_i(T,B+h,c_j,\tau)\le\nu<\nu_i(T,B,c_j,\tau);\\
\beta\sum_kP_{j,k}g_{h}^\nu(T-1,B,c_k,\{\tau+1\}),\\
\mbox{if~}\nu_i(T,B,c_j,\tau)\le\nu.
\end{array}
\right.
\end{array}
\]
\eitem

Denote
\[
\underline{\nu}_h(T,B,c_j,\tau)\triangleq\min_{\{k:P_{j,k}>0\}}\{\underline{\nu}_h(T-1,(B-1)^+,c_k,\tau)\},
\]
and
\[
\bar{\nu}_h(T,B,c_j,\tau)\triangleq\max_{\{k:P_{j,k}>0\}}\{\bar{\nu}_h(T-1,B,c_k,\tau)\},
\]
where $a^+=\max\{0,a\}.$

Since Proposition~\ref{Prop:differenceOfValueFunction1} holds for ${T-1}$ by assumption, we have ${g_h^\nu(T,B,c_j,\tau)}$ is continuous in $\nu$. When ${\nu\in[\underline{\nu}_h(T,B,c_j,\tau),\bar{\nu}_h(T,B,c_j,\tau)]}$, $g_h^\nu(T,B,c_j,\tau)$ is piece-wise linear and ${\partial g_h^\nu(T,B,c_j,\tau)/\partial\nu\ge -h}$. Otherwise, $g_h^\nu(T,B,c_j,\tau)$ is constant.
\end{proof}

\begin{proof}[Proof of Proposition~\ref{Prop:concavityOfValueFunction1}]\\
Since $\nu_i(T,B,c_j,\tau)$ is nondecreasing in $B$ when ${B\ge T}$,
\[\hspace{-1em}
\begin{array}{l}
\mathrel{\phantom{=}}g_1^\nu(T,B+1,c_j,\tau)-g_1^\nu(T,B,c_j,\tau)\\
=V_i^\nu(T,B+2,c_j,\tau)+V_i^\nu(T,B,c_j,\tau)-\\
\mathrel{\phantom{=}}2V_i^\nu(T,B+1,c_j,\tau)\\
=\left\{
\begin{array}{l}
\beta\sum_kP_{j,k}g_1^\nu(T-1,B,c_k,\{\tau+1\})-\\
\beta\sum_kP_{j,k}g_1^\nu(T-1,B-1,c_k,\{\tau+1\}),\\
\mbox{if~}\nu\le\nu_i(T,B,c_j,\tau);\\
\nu-(1-c_j)+\beta\sum_kP_{j,k}g_1^\nu(T-1,B,c_k,\{\tau+1\}),\\
\mbox{if~}\nu_i(T,B,c_j,\tau)<\nu\le\nu_i(T,B+1,c_j,\tau);\\
1-c_j-\nu-\beta\sum_kP_{j,k}g_1^\nu(T-1,B,c_k,\{\tau+1\}),\\
\mbox{if~}\nu_i(T,B+1,c_j,\tau)<\nu\le\nu_i(T,B+2,c_j,\tau);\\
\beta\sum_kP_{j,k}g_1^\nu(T-1,B+1,c_k,\{\tau+1\})-\\
\beta\sum_kP_{j,k}g_1^\nu(T-1,B,c_k,\{\tau+1\}),\\
\mbox{if~}\nu_i(T,B+2,c_j,\tau)\le\nu;\\
\end{array}
\right.\\
\le0.
\end{array}
\]
\normalsize
The first and forth cases are non-positive because ${V_i^\nu(T-1,B,c_j,\tau)}$ is concave when ${B\ge T-1}$ according to the assumption. The second and third cases are non-positive because of the expressions of $\nu_i(T,B+1,c_j,\tau)$. Thus $V_i^\nu(T,B,c_j,\tau)$ is concave in $B$ when $B\ge T$.
\end{proof}

\subsection{Proof of Optimality of Whittle's Index with $M=N$}\label{proof:optimalityOfWhittleIndexPolicy}
In this appendix, we prove
that the Whittle's index policy optimally solves the MAB problem defined in (\ref{eqn:MAB}), which is equivalent to the  MDP problem formulated in (\ref{eqn:originPro}).

First, we claim that the Whittle's index policy optimally solves the single charger problem with dynamic price and no constraint, as defined in (\ref{eqn:noConstraint}). The extended state $\tilde{s}=(T,B,c,\tau)$ includes the charging cost and period index. The Bellman equation of the single charger problem is given by:
\begin{equation}\label{eqn:singleArmMAB}
V_i(\tilde{s})=\max\{R_0(\tilde{s})+\beta(\mathcal{L}_{0}V_i)(\tilde{s}), R_1(\tilde{s})+\beta(\mathcal{L}_{1}V_i)(\tilde{s})\},
\end{equation}
 where action $a=1$ means to activate the charger and $a=0$ means to leave it passive.

The Whittle's index is defined by introducing a  {\it $\nu$-subsidy problem},
which is a modified version of the single arm problem defined in (\ref{eqn:singleArmMAB}).
In the $\nu$-subsidy problem, whenever the passive action is taken, the scheduler receives an extra reward $\nu$  \cite{Whittle:1988JAP}. The single charger problem defined in (\ref{eqn:singleArmMAB}) is simply the case when the subsidy $\nu=0$.

The Bellman equation for the $\nu$-subsidy problem is given by
\beq\label{eqn:nu-problem}
V_i^\nu(\tilde{s})=\max\{R_0(\tilde{s})+\nu+\beta(\mathcal{L}_{0}V_i^\nu)(\tilde{s}), R_1(\tilde{s})+\beta(\mathcal{L}_{1}V_i^\nu)(\tilde{s})\},
\eeq
where $V_i^\nu$ is the value function for the $\nu$-subsidy problem.

Now define a Whittle's index policy $\pi_1$ for a single charger (either regular or dummy charger) $\nu$-subsidy problem as to activates the charger if and only if $\nu_i(\tilde{s})>\nu$. Thus we have the following lemma.

\begin{lem}\label{col:optimalityOfWhittleindex}
The Whittle's index policy $\pi_1$ is optimal for the single charger $\nu$-subsidy problem defined in (\ref{eqn:nu-problem}). In particular, when ${\nu=0}$, $\pi_1$ is optimal for the single charger problem defined in (\ref{eqn:singleArmMAB}).
\end{lem}

\begin{proof}
We have shown in Appendix \ref{proof:indexability} that
the Whittle's index defined in Definition \ref{def:W} exists,
and therefore the Whittle's index policy $\pi_1$ is well defined.
By Definition \ref{def:W}, for any state $\tilde{s}$ such that ${\nu_i(\tilde{s})>\nu}$, the first term in (\ref{eqn:nu-problem}) is strictly smaller than the second term. The Whittle's index policy $\pi_1$ activates the charger  and obtains the second term as an expected reward which satisfies the Bellman equation in this case.

For ${\nu=\nu_i(\tilde{s})}$, the first term is greater or equal to the second term in the Bellman equation by Definition \ref{def:W}. The indexability result proved in Appendix \ref{proof:indexability} guarantees that the passive set grows monotonously which implies that this inequality is true for any ${\nu\ge\nu(\tilde{s})}$. Thus, for any state $\tilde{s}$ such that ${\nu(\tilde{s})\le\nu}$, the Whittle's index policy $\pi_1$ leaves the charger passive and obtains the first term as the expected reward, satisfying the Bellman equation.

Thus, $\pi_1$ satisfies the Bellman equation (\ref{eqn:nu-problem}) and is therefore optimal for the single charger $\nu$-subsidy problem. In particular, when ${\nu=0}$, $\pi_1$ is optimal for the single charger problem and satisfies the Bellman equation (\ref{eqn:singleArmMAB}).
\end{proof}

Now we consider the problem (\ref{eqn:MAB}) with ${M=N=1}$:
 we have a regular charger and a dummy charger, and at each time, we are required to activate exact one charger.
For this constrained two-arm problem, the state is defined as ${(\tilde{s},\textbf{0})=(T,B,c,\tau,0,0)}$, where $\tilde{s}$ is the extended state of the regular charger and $\textbf{0}=(0,0)$ the state of the dummy charger. The action ${a'=1}$ means to activate the regular charger, and ${a'=0}$ represents activating the dummy charger.

The state of the dummy charger will always be $\textbf{0}$. The dummy charger yields no reward regardless of the taken action. Thus the state transition of two-arm problem
is equivalent to the state transition in problem (\ref{eqn:singleArmMAB}), i.e., ${P((\tilde{s},\textbf{0}),(\tilde{s}',\textbf{0})|a')=P(\tilde{s},\tilde{s}'|a)}$.
The rewards of the two-arm problem can be presented by the rewards of
 the single charger problem  in (\ref{eqn:singleArmMAB}):
\[\begin{array}{l}
R'_1(\tilde{s},\textbf{0})=R_1(\tilde{s}),\\
R'_0(\tilde{s},\textbf{0})=R_0(\tilde{s}).
\end{array}
\]

The Whittle's index policy for the two-arm problem (denoted by $\pi_2$) activates the regular charger when ${(\nu(\tilde{s})>\nu(\textbf{0})=0)}$, and activates the dummy charger (leaving the regular charger passive) otherwise.

 When $\pi_1$ faces state $\tilde{s}$ and $\pi_2$ faces state $(\tilde{s},\textbf{0})$ for the same realization $\tilde{s}$, the actions of two policies are the same.  $\pi_2$ will activate the regular charger in the two-arm problem if and only if $\pi_1$ activates the charger in the single charger problem, and vice versa. Since the reward, transition and the action of these two policies are the same, the value functions will be the same. Denoting the value function of $\pi_1$ and $\pi_2$ by $V_{\pi_1}(\tilde{s})$ and $H_{\pi_2}(\tilde{s},\textbf{0})$, we have $H_{\pi_2}(\tilde{s},\textbf{0})=V_{\pi_1}(\tilde{s})$. Since  $V_{\pi_1}(\tilde{s})$ satisfies the Bellman equation (\ref{eqn:singleArmMAB}), we have
\[\begin{array}{l}
\mathrel{\phantom{=}}H_{\pi_2}(\tilde{s},\textbf{0})\\
=\max\{R_0(\tilde{s})+\beta(\mathcal{L}_{0}H_{\pi_2})(\tilde{s},\textbf{0}), R_1(\tilde{s})+\beta(\mathcal{L}_{1}H_{\pi_2})(\tilde{s},\textbf{0})\}\\
=\max\{R'_0(\tilde{s},\textbf{0})+\beta(\mathcal{L}_{0}H_{\pi_2})(\tilde{s},\textbf{0}),\\ \mathrel{\phantom{=\max\{}}R'_1(\tilde{s},\textbf{0})+\beta(\mathcal{L}_{1}H_{\pi_2})(\tilde{s},\textbf{0})\},
\end{array}\]
which is in fact the Bellman equation for the constrained two-arm problem.
The Whittle's index policy satisfies the Bellman equation for the two-arm problem and is therefore optimal.

Finally, we argue that the Whittle's index policy is optimal for the multi-arm problem defined in (\ref{eqn:MAB}).
We have $N$ regular chargers and $N$ dummy chargers. At each time, we activate exact $N$ chargers. We can pair each regular charger with¡¡a dummy charger and implement the Whittle's index policy for each pair. The action of each regular charger is decoupled, and the total reward is simply the sum of reward from all the $N$ regular chargers. The Whittle's index policy optimally optimal solves the problem of each pair, and is therefore optimal for the original problem in (\ref{eqn:MAB}).
We note, however, that  the above argument no longer holds when ${M<N}$, because the  problem defined in (\ref{eqn:MAB}) cannot be decoupled
into $N$ single (regular) charger problems in this case.

\old{
An alternative way is to consider the relationship between the single arm problem defined in (\ref{eqn:noConstraint}) and the problem in (\ref{eqn:originPro}), which is equivalent to the problem (\ref{eqn:MAB}). In problem (\ref{eqn:originPro}), when ${M=N}$, there is in fact no constraint (since ${\sum_ia_i\le N}$ always holds) and all chargers are decoupled. Instead of one big problem with $N$ chargers we have $N$ small problem of one charger. The reward is simply the sum of reward from individual chargers. Since Whittle's index policy optimizes the  individual charger problem, it optimally solves the whole problem.}

\subsection{Proof of Closed-form of Whittle's Index}\label{proof:closed-form}
\begin{proof}
Since the cost $c_0$ is constant, we will omit the cost in the state of chargers for simplicity.

For dummy chargers, there is no EV arrival, and only the charging cost evolves.
The Bellman equation of the $\nu$-subsidy problem is given by
\[\begin{array}{l}
V_i^\nu(0,0,\tau)=\max\{\nu+\beta\sum_kP_{j,k}V^\nu_i(0,0,\{\tau+1\}),\\[2pt]
\mathrel{\phantom{V_i^\nu(0,0,\tau)=\max\{}}\beta\sum_kP_{j,k}V^\nu_i(0,0,\{\tau+1\})\}.
\end{array}\]
When ${\nu<0}$, it is optimal to activate the dummy charger. Otherwise, passive action is optimal. So a dummy charger is indexable and its Whittle's index is $\nu_i(0,0,\tau)=0$.

For regular chargers, we showed in Appendix B2 that ${\nu_i(1,0,\tau)=0}$ and ${\nu_i(1,B,\tau)=1-c_0+F(B)-F(B-1)}$ when $B\ge1$. We will show the index closed-form for the case of ${T\ge2}$ using induction.

\subsubsection {When $T=2$} The discussion is divided into following two conditions.

If $B=1$,
 \[
 \begin{array}{l}
V_i^\nu(2,1,\tau)=\max\{\nu+\beta V_i^\nu(1,1,\{\tau+1\}),\\
 \mathrel{\phantom{V_i^\nu(2,1,\tau)=\max\{}}1-c_0+\beta V_i^\nu(1,0,\{\tau+1\})\}.
 \end{array}
 \]
 The difference between active and passive actions
\[
\begin{array}{l}
\mathrel{\phantom{=}}f^\nu(2,1,\tau)\\
=\nu-(1-c_0)+\beta g_1^\nu(1,0,\{\tau+1\})]\\
=\left\{
\begin{array}{ll}
\nu-(1-\beta)(1-c_0),& \mbox{if~}\nu<0;\\
(1-\beta)[\nu-(1-c_0)],& \mbox{if~}0\le\nu<1-c_0+F(1);\\
\nu-(1-c_0)-\beta F(1),& \mbox{if~}1-c_0+F(1)\le\nu;\\
\end{array}
\right.
\end{array}
\]
\normalsize
equals to $0$ when ${\nu=1-c_0}$. Thus ${\nu_i(2,1,\tau)=1-c_0}$.

If $B\ge2$, the Bellman equation is stated as follows.
\[\hspace{-1em}
\begin{array}{l}
V_i^\nu(2,B,\tau)=\max\{\nu+\beta V_i^\nu(1,B,\{\tau+1\}),\\
\mathrel{\phantom{V_i^\nu(2,B,\tau)=\max\{}}1-c_0+\beta V_i^\nu(1,B-1,\{\tau+1\})\}.
\end{array}
\]
Denote ${\Delta F(B)=F(B)-F(B-1)}$. The difference between active and passive actions
\[
\begin{array}{l}
\mathrel{\phantom{=}}f^\nu(2,B,\tau)\\
=\nu-(1-c_0)+\beta g_1^\nu(1,B-1,\{\tau+1\})\\
=\left\{
\begin{array}{l}
\nu-(1-c_0)-\beta\Delta F(B-1),\\
\mbox{if~}\nu<1-c_0+\Delta F(B-1);\\
(1-\beta)[\nu-(1-c_0)],\\
\mbox{if~}1-c_0+\Delta F(B-1)\le\nu<1-c_0+\Delta F(B);\\
\nu-(1-c_0)+\beta\Delta F(B),\\
\mbox{if~}1-c_0+\Delta F(B)\le\nu;\\
\end{array}
\right.
\end{array}
\]
equals to $0$ when ${\nu=1-c_0+\beta[F(B-1)-F(B-2)]}$. Thus ${\nu_i(2,B,\tau)=1-c_0+\beta[F(B-1)-F(B-2)]}$.

\subsubsection {When $T>2$}
Assume Equation (\ref{eqn:closedForm}) holds for ${T-1}$, consider the case for $T$.

If $B=1$,
\[\begin{array}{l}
V_i^\nu(T,B,\tau)=\max\{\nu+\beta V_i^\nu(T-1,1,\{\tau+1\}),\\
\mathrel{\phantom{V_i^\nu(T,B,\tau)=\max\{}}1-c_0+\beta V_i^\nu(T-1,0,\{\tau+1\})\}.
\end{array}
\]
The difference between actions is
\[\hspace{-1em}
\begin{array}{l}
\mathrel{\phantom{=}}f^\nu(T,1,\tau)\\
=\nu-(1-c_0)+\beta g_1^\nu(T-1,0,\{\tau+1\})\\
=\left\{
\begin{array}{ll}
\nu-(1-\beta)(1-c_0),& \mbox{if~}\nu<0;\\
(1-\beta)[\nu-(1-c_0)],& \mbox{if~}0\le\nu<1-c_0;\\
\nu-(1-c_0)+\\
\beta^2g_1^\nu(T-2,0,\{\tau+2\})& \mbox{if~}1-c_0\le\nu.\\
\end{array}
\right.
\end{array}
\]
The last case can be rewritten as
\[\hspace{-1em}
\begin{array}{l}(1-\beta)[\nu-(1-c_0)]+\beta(\nu-(1-c_0))+\\
\beta^2(V_i^\nu(T-2,1,\{\tau+2\},\nu)-V_i^\nu(T-2,0,\{\tau+2\})),
\end{array}\]
which equals to $0$ when ${\nu=1-c_0}$ since ${\nu_i(T-1,1,\tau)=1-c_0}$ by assumption.  Thus ${\nu_i(T,1,\tau)=1-c_0}$.

If $2\le B\le T-2$, the difference between actions is stated as follows.
\[
\begin{array}{l}
\mathrel{\phantom{=}}f^\nu(T,B,\tau)\\
=\nu-(1-c_0)+\beta g_1^\nu(T-1,B-1,\{\tau+1\})\\
=\left\{
\begin{array}{ll}
\nu-(1-c_0)+\\
\beta^2g_1^\nu(T-2,B-2,\{\tau+1\})&\mbox{if~}\nu<1-c_0;\\
\nu-(1-c_0)+\\
\beta^2g_1^\nu(T-2,B-1,\{\tau+1\})&\mbox{if~}1-c_0\le\nu.
\end{array}
\right.
\end{array}
\]
The latter case equals to $0$ when ${\nu=1-c_0}$ since ${\nu_i(T-1,B,\tau)=1-c_0}$ when ${2\le B\le T-2}$ by assumption. Thus ${\nu_i(T,B,\tau)=1-c_0}$ when ${2\le B\le T-2}$.

If $B=T-1$,
\[\hspace{-1em}\begin{array}{l}
\mathrel{\phantom{=}}f^\nu(T,B,\tau)\\
=\nu-(1-c_0)+\beta g_1^\nu(T-1,B-1,\{\tau+1\})\\
=\left\{
\begin{array}{l}
\nu-(1-c_0)+\beta^2g_1^\nu(T-2,B-2,\{\tau+2\}),\\
\mbox{if~}\nu<1-c_0;\\
(1-\beta)[\nu-(1-c_0)],\\
\mbox{if~}1-c_0\le\nu<1-c_0+\beta^{T-2}F(1);\\
\nu-(1-c_0)+\beta^2g_1^\nu(T-2,B-1,\{\tau+2\}),\\
\mbox{if~}1-c_0+\beta^{T-2}F(1)\le\nu;
\end{array}
\right.
\end{array}
\]
equals to $0$ when $\nu=1-c_0$. So ${\nu_i(T,B,\tau)=1-c_0}$ when $B=T-1$.

If $B\ge T$,
\begin{equation}\label{eqn:differenceOfValueFunction}
\begin{array}{l}
\mathrel{\phantom{=}}f^\nu(T,B,\tau)\\
=\nu-(1-c_0)+\beta g_1^\nu(T-1,B-1,\{\tau+1\})\\
=\left\{
\begin{array}{l}
\nu-(1-c_0)+\beta^2g_1^\nu(T-2,B-2,\{\tau+1\})\\
\mbox{if~}\nu<1-c_0+\beta^{T-2}\Delta F(B-T+1);\\
(1-\beta)[\nu-(1-c_0)],\\
\mbox{if~}1-c_0+\beta^{T-2}\Delta F(B-T+1)\\
\mathrel{\phantom{\mbox{if~}}}\le\nu<1-c_0+\beta^{T-2}\Delta F(B-T+2);\\
\nu-(1-c_0)+\beta^2g_1^\nu(T-2,B-1,\{\tau+1\}),\\
\mbox{if~}1-c_0+\beta^{T-2}\Delta F(B-T+2)\le\nu.\\
\end{array}
\right.
\end{array}
\end{equation}
When $\nu<1-c_0+\beta^{T-2}[F(B-T+1)-F(B-T)]$,
\[
\begin{array}{l}
\nu\le\nu_i(T-1-T',B-1-T',\tau)\\
\mathrel{\phantom{\nu}}\le\nu_i(T-1-T',B-T',\tau)
\end{array}
\]
for all ${0\le T'\le T-1}$. Thus in the first case of (\ref{eqn:differenceOfValueFunction}),
\[\hspace{-1em}
\begin{array}{l}
\mathrel{\phantom{=}}\beta^2g_1^\nu(T-2,B-2,\{\tau+1\})\\
=\beta^3 g_1^\nu(T-3,B-3,\{\tau+2\})\\
=\cdots\\
=\beta^{T-1} g_1^\nu(1,B-T+1,\{\tau+T-2\})\\
=\beta^{T-1}[-F(B-T+1)+F(B-T)]
\end{array}
\]

 So when ${\nu=1-c_0+\beta^{T-1}[F(B-T+1)-F(B-T)]}$, the first case in equation (\ref{eqn:differenceOfValueFunction}) equals to $0$ . Thus when $B\ge T$, the closed-form of index is stated as: \[\nu_i(T,B,\tau)=1-c_0+\beta^{T-1}[F(B-T+1)-F(B-T)].\]
\end{proof}

%
%

{
\bibliographystyle{ieeetran}
\bibliography{Bibs/Journal,Bibs/Conf,Bibs/Book,Bibs/Misc}
}


\end{document}